\documentclass[12pt]{amsart}
\usepackage{amssymb, euscript}
\usepackage[arrow, matrix]{xy}  

\newtheorem{theorem}[subsubsection]{Theorem}
\newtheorem{propos}[subsubsection]{Proposition}

\newtheorem{coroll}[subsubsection]{Corollary}

\theoremstyle{definition}

\def\g{{\mathfrak{g}}}                     
\def\h{{\mathfrak{h}}}
\def\t{{\mathfrak{t}}}

\def\z{{\mathfrak{z}}}
\def\Q{{\mathfrak{Q}}}
\def\A{{\mathfrak{A}}}
\def\k{{\mathfrak{k}}}
\begin{document}

\title{On the Generalized Futaki Invariant}

\author{Mirroslav Yotov}
\address{I.H.E.S.\\ Le Bois-Marie\\35, Route de Chartres\\F-91440 Bures-sur-Yvette, France}
\email{yotov@ihes.fr}


\begin{abstract}
We study the algebraic properties of the generalized Futaki invariant
of an almost Fano variety and prove that it is in fact a pushforward
to a point of an appropriate equivariant Chow class of the
variety. This allows us to use Bott-type formulae for calculating the
invariant. We show this use on some examples.
\end{abstract}

\maketitle

\section{Introduction}

Let $M$ be a compact complex manifold with positive first Chern
class. Such manifolds are called Fano manifolds. The manifold $M$\, is
called Einstein-K\"ahler if there exists a K\"ahler metric $g$\, on
it, such that the K\"ahler form $\omega _g$\, and the Ricci form $\rho
_g$\, of $g$\, satisfy the equation : $\omega _g = \rho _g$\,. In 1957
Matsushima \cite{Mats} has proved that the existence of
Einstein-K\"ahler metrics 
on $M$\, implies that the algebra ${\A}(M)$\, of holomorphic
vector fields on $M$\, is reductive. In 1983 A. Futaki \cite{Fut1}
introduced a 
character $F$\, of ${\A}(M)$\, which depends only on the complex
structure of $M$\,, and vanishes identically if $M$\, is
Einstein-K\"ahler. Futaki also showed an example of a Fano manifold
with reductive algebra ${\A}(M)$\, and with nontrivial character
$F$\,. $F$\, is called the Futaki invariant of the Fano manifold
$M$\,.

In spite of its analytic definition, $F$\, has very algebraic in
nature properties. For example, Mabuchi \cite{Mab1} has proved that
$F$ \,vanishes 
on all nilpotent elements of ${\A}(M)$\,. Secondly, Futaki and
Morita \cite{F-M} have shown that $F$\, can be defined by using appropriate
polynomials invariant under the action of ${\A}(M)$\, (we call
these Futaki-Morita polynomials), and have proved that there is a
relation between $F$\, and the equivariant cohomologies of $M$\,.

In 1992 Ding and Tian \cite{D-T} defined a generalization of the
Futaki invariant 
for almost Fano varieties $X$\,(i.e., normal complete ${\bf Q}$-Gorenstein
varieties with ample degrees of their anticanonical sheaf). They
proved that, if one considers special degenerations of a Fano manifold
$M$\,(such degenerations have almost Fano central fibres), then the
existence of Einstein-K\"ahler metrics on $M$\, depends on the
behaviour of the generalized Futaki invariants of the central fibres
of these degenerations. This generalized Futaki invariant proved to be
a very efficient tool in investigating Fano manifolds. By using a
refined version of Ding-Tian's theorem, Tian \cite{Tian} found a
counterexample to 
the long-standing hypothesis saying that each Fano manifold with a
discrete group of biholomorphisms is Einstein-K\"ahler. The
significance of the generalized Futaki invariant stems also from the
fact that it is connected, via the notion of weakly K-stability of a
Fano manifold introduced by Tian, with the Chow-Mumford stability
w.r.t. the very ample degrees of the anticanonical sheaf of $X$\,
\cite{Tian}.
In this paper we study the properties of the generalized Futaki
invariant in the spirit of Futaki-Morita and Mabuchi. Establishing the
connection of $F$\, with the equivariant cohomologies of $X$\,, we
found it very convenient to express $F$\, in terms of Edidin-Graham's
equivariant Chow cohomology groups. In particular. we use the
definition of equivariant Chern classes as given in \cite{E-G1}, and
avoid using of invariant connections to define them. When the group
acting is compact, these coincide with the equivariant characteristic
classes as defined by Berline and Vergne \cite{B-V} . As a result we
prove that $F$\,can 
be regarded as a pushforward to a point of an equivariant Chow
cohomology class on $X$\,. The benefit from this result is
two-fold. First of all, it gives a natural way of defining a Futaki
invariant for the varieties other than normal almost Fano ones - this
is useful when studying the degenerations of Fano
manifolds. Secondly, it allows us to use the machinery of Bott-type
residue formulas for calculating $F$\,. In this way we recover not
only the known (Bott-type formulas) for $F$\, but also the one for
complete intersections due to Zhiqin Lu \cite{Lu}.

$ $\par

We now explain the content of the paper.

In Section 2 we recall some basic properties of the Edidin-Graham's
equivariant Chow groups of an algebraic scheme in the extent needed
for our purposes. We emphasize here  the behaviour of these groups
when passing to a subgroup acting on the scheme. We find a way of
attaching to each cohomology class $\alpha $\,on $X$\, a map
$f_{\alpha }$\, 
from the Lie algebra of the group of biholomorphisms of $X$\, to the
complex numbers, and show that $f_{\alpha }$\, is actually a character
of this Lie algebra for some classes $\alpha $\,.

In Section 3 we recall the definition of the generalized Futaki
invariant $F$\, of an almost Fano variety $X$\,, show that it is
insensitive to the singularities of $X$\, and prove two useful
integral representations. One of these, combined with an argument of
Mabuchi \cite{Mab1}, helps one to prove that $F$\, vanishes on the
nilpotent 
elements of $Lie\,Aut(X)$\,. The other shows that $F$\, can be defined
by using an appropriate Futaki-Morita polynomial.

Section 4 is devoted to a study of the Futaki-Morita polynomials for
non-smooth algebraic varieties. In this section we find a necessary
and sufficient condition for a Futaki-Morita polynomial to be
represented by some $f_{\alpha }$\,. In particular we find such
$\alpha $\, for the generalized Futaki invariant.

In section 5 we recall the essence of the Bott-type residue formulas,
as dive-loped by Edidin-Graham, and give some examples of calculating
$F$\,. We first, without appeal to Bott formulas, prove the Lu's
formula for $F$\, when $X$\, is a complete intersection in a
projective space. Then we use the Bott formula for calculating for
some degenerations of the blow-up of the three dimensional projective
space in a twisted cubic curve.

$ $\par

{\bf Acknowlegements.} This paper was written while the author was an
EPDI fellow and enjoyed the hospitality of I.H.E.S. and Fourier
Institute in Grenoble. He is grateful to both institutions for giving
him the opportunity to work in a very inspiring atmosphere! The author
is grateful to M.Brion and J.-P.Demailly for the fruitful discussions
he has had with them! Special thanks are due to T.Pantev who helped
the author to make the preliminary variant of this paper readable as
well as for making him aware of the preprint by Zhiqin Lu \cite{Lu}. 

$ $\par

\section{Basic facts from the theory of the Equivariant Chow groups}

The purpose of this paper is to establish the algebraic nature (in terms of 
equivariant Chow groups) of both - the classical Futaki invariant \cite{Fut1} 
and the generalized Futaki invariant \cite{D-T}, and to give some examples of 
how to calculate them. In this section we present some definitions and recall 
some properties of the equivariant Chow groups to the extent and in form we 
will need them in the text. The aim is to define, for each element $\alpha $ 
of an operational equivariant Chow group of a normal algebraic variety, a 
polynomial $f_{\alpha}$ on the Lie algebra of the group of automorphisms of 
that variety. In section 3 we show that both the classical and the generalized
 Futaki invariants are special cases of this construction.

The main sources on the theory of Chow groups for us are \cite{Ful}, 
\cite{E-G1}, and \cite{Bri1}. We work over the field of complex numbers 
${\bf C}$.

\subsection{Equivariant Chow Groups}

In our research we will use the definition of equivariant Chow groups
due to  Edidin and Graham \cite{E-G1}. By an algebraic group we always
mean a linear algebraic group. 

Let $X$ be a scheme (over ${\bf C}$), \, $\dim X = n$, and let $G$\,
be an algebraic group acting on $X$\,. Suppose $V$\, is a linear
representation of $G$\, such that there is an open subset $U \subset
V$\, on which $G$\, acts  freely, and let $U \rightarrow
 U/G$\, be the
principal quotient bundle. It is shown in \cite{E-G1} that for any
(linear) algebraic group $G$\, and for any $q$\, there exists a linear
representation $V$\, of $G$\, such that $codim_V(V-U) > n-q$\, and the
quotient space $U/G$\, is a scheme. The mixed space $X_G :=
X\times_GU$\,\, exists as an algebraic space, but if $X$\, is normal,
the case we are mainly interested in, $X_G$\, is in fact a scheme. We
define the $q$-th $G$-equivariant Chow group of $X$\, by setting
$$A^G_q(X) := A_{q+\dim V-\dim G}(X_G),$$
where on the right-hand-side we use the usual Chow group. As is shown
in \cite{E-G1}  the definition above does not depend on the particular
choice of the representation. Obviously, $A_q^G(X) = 0$\,\, for all $q
> \dim X$. 

Without going into the theory of the equivariant Chow groups ( for
which we refer to \cite{E-G1} and \cite{Bri1}) we recall here the
properties which govern their behaviour when passing to a closed subgroup. 

$ $\par

1. Let $H < G$\, be a closed subgroup of $G$\,, and let $V$\, be a 
representation of $G$\, which defines the group $A_q^G(X)$\, as above. Hence,  $V$\, defines also the group $A_q^H(X)$\,. There is a natural smooth morphism
$$f_{H,G} : X\times_HU \longrightarrow X\times_GU$$
with fibre $G/H$\,. The flat pull back morphism for the usual Chow groups
$$f_{H,G}^* : A_*(X\times_GU) \longrightarrow A_{*+\dim G/H}(X\times_HU).$$   
can be interpreted as a morphism
$$\phi_*^{H,G} : A_*^G(X) \longrightarrow A_*^H(X)$$
of equivariant Chow groups.

In particular, if $G/H$\, is isomorphic to an affine space, then $\phi_*^{H,G}$\, will be an isomorphism (\cite{Gill}) due to the homotopy invariance of the usual Chow groups. As a consequence we get that $\phi_*^{H,G}$\, is an isomorphism for any maximal reductive closed subgroup $H$\, of $G$. (We call such subgroups Levi subgroups of $G$).

$ $\par

{\bf Remark.}  In the notations above, $X$\, defines elements $[X]_G$\, and 
$[X]_H$\, in $A_n^G(X)$\, and $A_n^H(X)$\, respectively. We have that 
$\phi_*^{H,G}([X]_G) = [X]_H$.

$ $\par

2. (Proper push-forward)  For any proper $G$-scheme $X$\,, the canonical map
$$\pi_X : X \longrightarrow pt,$$
viewed naturally as a $G$-equivariant map, gives rise to a morphism
$$ \pi^G_{X*} : A_*^G(X) \longrightarrow A_*^G(pt).$$
This map can be seen as an integration along the fibres of the fibration
$$X\times_GU \longrightarrow pt\times_GU = U/G$$
with fibre $X$.

The commutative diagram of maps
$$\xymatrix{
X\times_HU \ar[r]^{f_{H,G}} \ar[d] & X\times_GU \ar[d]\\
pt\times_HU \ar[r] & pt\times_GU}$$
gives rise to the commutative diagram
$$\xymatrix{
A_*^G(X) \ar[r]^{\phi_*^{H,G}} \ar[d]_{\pi^G_{X*}}
& A_*^H(X) \ar[d]^{\pi ^H_{X*}} \\
A_*^G(pt) \ar[r]^{\phi_*^{H,G}} & A_*^H(pt) ,}$$
i.e. $ \pi^H_{X*} \circ \phi_*^{H,G} = \phi_*^{H,G} \circ \pi^G_{X*}$.

$ $\par

3. The $q$-th $G$-equivariant Borel-Moore homology group is defined as
$$H^G_{BM,q}(X) := H_{BM,q+2\dim V-2\dim G}(X_G)$$
notations and assumptions being as above. There is a natural ``cycle'' map 
( cf. \cite{E-G1})
$$cl^G_X : A_q^G(X) \longrightarrow H^G_{BM,2q}(X),$$
compatible with the operations on the equivariant Chow groups. As in the case 
of the Chow groups there are morphisms 
$$\phi_*^{H,G} : H^G_{BM,q}(X) \longrightarrow H^H_{BM,q}(X),$$
such that the diagram
$$\xymatrix{
A_q^G(X) \ar[r]^{cl_X^G} \ar[d]_{\phi_*^{H,G}}
                         & H^G_{BM,2q}(X) \ar[d]^{\phi_*^{H,G}}\\
A_q^H(X) \ar[r]^{cl_X^H} & H^H_{BM,2q}(X)}$$
commutes, i.e. $\phi_*^{H,G}\circ cl_X^G = cl_X^H\circ \phi_*^{H,G}$.

$ $\par

{\bf Remark.} When $X$\, is smooth $H_{BM,q}^G(X) \cong H^{2n-q}(X_G) \cong 
H^{2n-q}_G(X)$\, the last being the $G$-equivariant cohomology group of $X$.

$ $\par

\subsection{Equivariant Operational Chow groups}

We will use the equivariant operational Chow groups as defined in \cite{E-G1}. For $i \geq 0$\, we set
$$ A^i_G(X) = \left\{ c(f) : A_*^G(Y) \rightarrow A_{*-i}^G(Y) \; \left| \; 
\begin{minipage}[c]{5cm}
$f : Y \to X$ is a $G$-equivariant map of normal varieties
\end{minipage}
\right. \right\}$$  
where the morphisms $c(f)$\, are compatible with the operations on the 
equivariant Chow groups (pull-back for l.c.i. morphisms, proper push-forwards, 
etc.). The composition of maps defines a graded ring structure on 
$A^*_G(X) := \bigoplus_{i\geq 0} A_G^i(X)$.

We list below some of the properties of the equivariant operational Chow 
groups that we will be using further.

$ $\par

1. If $V$\, is a linear representation of $G$\,, $U\subset V$\, is an open 
subset on which $G$\, acts freely, and codim$_V(V-U) > k$\,, then
$$A^k_G(X) = A^k(X_G).$$

Let $E \rightarrow X$\, be a $G$-equivariant vector bundle. Then $E_G := 
U\times_GE$\, is a vector bundle over $X_G$\,. The equivariant Chern class 
$c_i^G(E)$\,of $E$\, is defined as the operational class $c_i(E_G)$.

$ $\par 

2. (Cap-product)  For any $\alpha \in A_G^q(X)$\,, the identity map $id : X \rightarrow X$\, defines an operation $\alpha(id) : A_p^G(X) \rightarrow A_{p-q}^G(X)$\, determining in this way the cap-product
$$\cap \, : \, A_G^q(X) \otimes A_p^G(X) \longrightarrow A_{p-q}^G(X)$$
$$\alpha \otimes a \, \longmapsto \, \alpha \cap a \, := \, \alpha (id)(a).$$
This cap-product makes $A_*^G(X)$\,  into a $A^*_G(X)$\,-module.

$ $\par

3. The canonical map $\pi _X^G : X \rightarrow pt$\, defines a degree $0$ 
morphism of graded rings
$$\pi_X^{G*} : A^*_G(pt) \longrightarrow A^*_G(X).$$
Under this morphism $A_*^G(X)$\, becomes a $A^*_G(pt)$- module.

$ $\par 

4. Let $H < G$\, be a closed subgroup. The natural morphism, defined as above,
$$U\times_HX \longrightarrow U\times_GX$$
defines a morphism
$$A^*(X_G) \longrightarrow A^*(X_H)$$
which, when codim$_V(V-U)$\, is appropriate, gives rise to a morphism
$$\phi_{H,G}^* : A_G^*(X) \longrightarrow A_H^*(X).$$
Further, the morphisms $\phi_{H,G}^*$\, and $\phi^{H,G}_*$\, are compatible 
with the cap-product morphism - the diagram
$$\xymatrix{
A_G^k(X)\otimes A_m^G(X) \ar[r]^{\cap} \ar[d]_{\phi^*_{H,G}\otimes \phi_*^{H,G}}
                                       & A_{m-k}^G(X) \ar[d]^{\phi_*^{H,G}} \\
A_H^k(X)\otimes A_m^H(X) \ar[r]^{\cap} & A_{m-k}^H(X)  }$$
commutes.

$ $\par

5. When $X$\, is smooth, the mapping
$$A_G^q(X)\, \ni \, c \longmapsto c\,\cap [X]_G \, \in \, A^G_{n-q}(X)$$
is an isomorphism, i.e. $A^*_G(X) \cong A^G_{n-*}(X)$. Due to this
there is an induced  ``cycle class map'' 
$$cl^X_G : A^*_G(X) \longrightarrow H^{2*}_G(X)$$
which commutes with the maps $\phi _{H,G}^*$\,:
$$\phi _{H,G}^*\circ cl_G^X = cl^X_H\circ \phi ^*_{H,G}.$$

$ $\par 

6. When $X$\, is a complete $G$- scheme there is an integration map
$$\int_X^G :=  \pi_*^G \circ ( \bullet  \cap [X]_G) : \, A_G^q(X) \longrightarrow A_{n-q}^G(pt).$$
For any closed subgroup $H < G$\,, we have
$$\int_X^H \circ \, \phi^*_{H,G} = \phi_*^{H,G} \circ \int_X^G.$$
In what follows we consider only complete schemes $X$\,, such that the
map $\int_X^G$\, is always defined.

$ $\par

7. (Structure Theorem) Suppose $G$\, is a reductive algebraic group,
$T < G$\, a maximal torus, and $ W:= N(T)/T$\, -the Weyl group of
$(G,T)$\,. Let $\chi (T)$\,  be the group of characters of $T$\,, and
denote by $S_{{\bf Z}}(T) = {\bf Sym}[\chi(T)]$\, and $S(T) = S_{{\bf
 Z}}(T)\otimes_{{\bf Z}}{\bf Q}$\, the corresponding symmetric algebras.
The Weyl group $W$\, acts naturally on  $A_*^T(X)$\, for any
$G$-scheme $X$\, (cf. \cite {E-G1} and \cite {Bri1}). We have the
following structure theorem

\begin{theorem}\label{struct}
The graded ring $A^*_T(pt)$\, is isomorphic to $S_{{\bf Z}}(T)$\,. The map
$$\phi^*_{T,G} : A_G^*(pt) : \longrightarrow A_T^*(pt)$$
is an injection over ${\bf Q}$\,, and identifies $A^*_G(pt)_{{\bf Q}} := A^*_G(pt)\otimes_{{\bf Z}}{\bf Q}$\, with the $W$-invariant subring $S(T)^W$\, of $S(T)$.

For any $G$-scheme $X$\,, the map $\phi_*^{T,G}$\, induces and isomorphism
$$A_*^G(X)_{{\bf Q}} \cong A_*^T(X)_{{\bf Q}}^W.$$

\end{theorem}  

$ $\par

\subsection{Definition and basic properties of the map $f_{\alpha}$}

Let $\xi \in \g $ be an element of the Lie algebra $\g$\, of $G$. Denote by $G(\xi)$\, the minimal algebraic subgroup of $G$\, containing $\exp(\xi)$\,. If $\xi = \xi _s + \xi _n$\, is the Jordan decomposition of $\xi$\,, 
then  $G(\xi) = T(\xi) \times G(\xi _n)$\,, where $T(\xi) = G(\xi _s)$\, is a 
subtorus of $G$\,, uniquely determined by $\xi$\,, and $G(\xi _n) = \{ \exp(t\xi _n) | t \in {\bf C}\}$\, is a unipotent subgroup of $G$. (cf. 
\cite{O-V}.)
\par

For any $\alpha \in A_G^q(X)$\, the element
$$\int^{G(\xi)}_X \phi^*_{G(\xi),G}(\alpha ) \, = \, \phi_*^{G(\xi),G}\circ \int_X^G \alpha$$
can be identified via $\phi_*^{T(\xi),G(\xi)}$\, with the element
$$\int_X^{T(\xi)} \phi^*_{T(\xi),G} (\alpha )\, \in \,
A_{n-q}^{T(\xi)}(pt) \,\cong \,A^{q-n}_{T(\xi)}(pt).$$
By the structure theorem above,
$$A^{T(\xi)}_{n-q}(pt) \, = \, {\bf Sym}[\chi (T(\xi))]_{q-n}$$
and $\int_X^{T(\xi)} \phi^*_{T(\xi),G}(\alpha )$\, is a polynomial of
degree $q-n$\, on $Lie\, T(\xi)$. (It is zero for $q < n$.)

Define
$$ f(\xi , \alpha) \, := \, ( \int _X^{T(\xi)} \phi^*_{T(\xi),G} (\alpha )) (\xi _s )\, \in \, {\bf C}.$$
The number $f(\xi ,\alpha )$\, depends only on $\xi$\, and $\alpha$\,, and for any $g\,\in \, G$\,
$$f(ad(g)\xi , \alpha ) \, = \, f(\xi ,\alpha ).$$

Another way of calculating $f(\xi ,\alpha )$\, is via maximal tori in $G$\, 
containing $T(\xi )$. \,Choose a maximal torus $T$ of $G$\, such that $T(\xi ) < T$.\, Then, $\xi _s \, \in \, Lie\, T(\xi ) \, \subset \, Lie\, T$\,\,, and
$$f(\xi ,\alpha ) = (\phi_*^{T(\xi ),T}\circ \phi _*^{T,G} \int_X^G \alpha )(\xi _s) $$
$$= (\phi _*^{T,G} \int_X^G \alpha )(\xi _s).$$ 

As is clear from their definition, the numbers $f(\xi ,\alpha )$\, depend 
linearly on the second argument. The dependence on the first argument is rather complicated in general. Nevertheless, there are some cases in which, for fixed second argument, $f_{\alpha } := f(\bullet ,\alpha )$\, is an 
$Ad(G)$-invariant polynomial on $\g$. One such case, which is closely related 
to the generalized Futaki invariant, is described below.

$ $\par

Let $T < G$\, be a maximal torus of $G$\,, let $H < G$\, be a Levi subgroup 
(i.e. a maximal closed reductive subgroup of $G$) containing $T$\,. Let $W := W(H, T)$\, be the corresponding Weyl group. Then,
$$G = H \ltimes Rad_uG,$$
$$H = Z(H).H',$$
$$T = (Z(H))_0\times T',$$
where $Rad_uG$ is the unipotent radical of $G$\,, $H' = (H,H)$\, is the commutator of $H$\,, $T' = T \cap H'$\, is the 
maximal torus of $H'$\,, corresponding to $T$\,, and $(Z(H))_0$\, is the 
identity component of the centre $Z(H)$\,. The Weyl group $W$\, can be 
naturally identified with $W(H',T')$\,.

We already know by the structure theorem that (over ${\bf Q}$\,) the map $\phi _*^{T,G}$\, is an embedding with image
$${\bf Sym}[\chi (T)]^W = {\bf Sym}[\chi (Z(H))]\otimes _{{\bf Q}}{\bf Sym}[\chi (T')]^W.$$
Suppose now that $\alpha $\, is such that
$$\phi^{T,G}_*\int^G_X \alpha \, \in \, {\bf Sym}[\chi (Z(H))]$$
and denote by $j : \g \rightarrow Lie\, Z(H)$\, the Lie algebra morphism 
corresponding to the canonical map $G \rightarrow H/H'$.\, Then $\alpha $\, 
defines an Ad$\,G$-invariant polynomial $f_{\alpha }$\, on $\g$ : for $ \xi  \, \in \, \g $
$$ f_{\alpha }(\xi ) := (\phi _*^{T,G} \int^G_X \alpha )(j(\xi )).$$ 

The algebra ${\bf Sym}[\chi (T')]^W$\, has no elements of degree one. Hence
$$({\bf Sym}[\chi (Z(H))]\otimes _{{\bf Q}}{\bf Sym}[\chi (T')]^W)_1 = ({\bf Sym}[\chi (Z(H))])_1,$$
and we have thus proved the following

\begin{propos}\label{prop}
In the above notations if $\alpha \, \in A_G^{n+1}(X)$\,, then
$f_{\alpha }$\, is a Lie algebras morphism 
$$f_{\alpha } : \g \longrightarrow {\bf C}$$
 and $ f_{\alpha }(\xi ) = f(\xi ,\alpha )$.\, This morphism vanishes
 on the nilpotent elements of $\g$.
  
\end{propos}

$ $\par

{\bf Remark} A more intrinsic proof of the proposition above, which I
owe to M.Brion, goes as follows. For any algebraic group $G$\, the
group $A^1_G(pt)$\, is naturally isomorphic to the Picard group
$Pic^G(pt)$\, of the $G$-equivariant line bundles over a point. The
latter is the same as the character group of $G$\,: each
$G$-equivariant line bundle
$L \rightarrow pt$\, corresponds to a homomorphism
$$\chi _L : G \rightarrow {\bf C}^*\,.$$
So, given $\alpha \in A^{n+1}_G(X)$\, its push-forward $\int^G_X
\alpha \in A^1_G(pt)\,$ can be identified with a character $\chi _L$\, of
$G$\,. The map $f_{\alpha }$\, is just the differential $d\chi _L$\,
of that character.

$ $\par

\section{The Generalized Futaki Invariant}

In this section we recall the definition of the generalized Futaki
 invariant  
 following \cite{D-T}, \cite{Tian}, and show that it has the 
properties of the map $f_{\alpha}$\, from Proposition \ref{prop}. In the case 
of the classical Futaki invariant this is shown by Futaki-Morita \cite{F-M}, 
and Mabuchi \cite{Mab1}. We follow very closely the approach in the cited 
papers with the appropriate changes needed when working with almost Fano 
varieties - the object of our interest. In the next section we'll show that in 
fact the properties from Proposition \ref{prop} are sufficient for the
 generalized Futaki invariant to be 
represented as $f_{\alpha }$\, for an appropriate $\alpha $\,.

$ $\par

Let $X$\, be a normal complete {\bf Q}-Gorenstein variety. Let $L
\rightarrow X$\, be the
line bundle, such that $L^*_{|X_{reg}} = \omega ^{\otimes k}_{X_{reg}}$\,, 
where $\omega_{X_{reg}}$\, is the dualizing sheaf of the regular part $X_{reg}$
\, of $X$\,. Suppose $L$ is ample. Such varieties $X$\, are called {\bf almost 
Fano varieties} (\cite{D-T}, \cite{Tian}).

A K\"ahler form $\omega $\, on $X_{reg}$\, is called {\bf admissible} ( on $X$
\,) if there exist an embedding of $X$\,, defined by some power of $L$\,,
$$\phi _{L^m} \, : \, X \, \hookrightarrow \, {\bf P}^N,$$
and a K\"ahler form $\tilde {\omega }$\, on ${\bf P}^N$\, representing 
$2\pi c_1({\bf P}^N)$\,, such that
$$\omega \, = \, {1\over km}.\phi ^*_{L^m} (\tilde \omega) \,\,\,\,\,\,\,\,\,\,\,\,\,\,\,\,\,\,\, {\text on} \,\,\,\,\,\,\,\,\,\,\,\,\,\, X_{reg}.$$

A vector field $\xi $ on $X_{reg}$\, is called {\bf admissible} (on $X$\,) if 
there exist an embedding as above, and a vector field $\tilde {\xi }$\, on 
${\bf P}^N$\,, such that $\tilde {\xi }$\, is tangent to $X$\, along $X_{reg}$
\,, and 
$$\xi \, = \, \tilde {\xi }_{|X_{reg}}.$$

Suppose $\omega $ is an admissible form on $X$\,. Then on $X_{reg}$\,
$$Ric(\omega ) \, - \, \omega \, = \, \sqrt{-1}\partial \bar {\partial }f_{
\omega },$$
where $Ric(\omega ) = -\sqrt{-1}\partial \bar {\partial }\log \omega ^n$\,, 
and $f_{\omega } \in C^{\infty}(X_{reg}, {\bf R})$. It is proved by Ding and 
Tian \cite{D-T} (see also \cite{Yoto}) that, given an admissible vector field 
$\xi $\,, for any admissible form $\omega $\, the number
$$F(\xi ) \, = \, {1\over (2\pi)^n }\,\int_{X_{reg}}\xi (f_{\omega }
) \omega ^n$$
is well defined, does not depend on the choice of $\omega $\,, and defines a 
character of the algebra of the admissible vector fields on $X$\,. When $X$\, 
is smooth, i.e. $X$\, is a Fano manifold, the map $F$\, coincides with the 
classical Futaki invariant \cite{Fut1}. In the general case of almost Fano 
varieties the map $F$\, is called {\bf the generalized Futaki
  invariant} of $X$. The 
significance of this invariant stems from the fact (\cite{D-T} and \cite{Tian})
 that it is closely related to the existence of Einstein-K\"ahler metrics on 
Fano manifolds.

$ $\par

Now we want to show that $F$\, has an appropriate integral representation, by 
using which one can prove that it vanishes on the nilpotent admissible vector 
fields, and can relate it to the Futaki-Morita polynomials from the next 
section.\par

Suppose $\omega = (1/km)\phi _{L^m}^*\tilde {\omega}$\,, for some
embedding $\phi_{L^m} : X \hookrightarrow {\bf P}^N$\,. Since the
admissible vector fields on $X$\, correspond to the elements of
$Lie\,G$\,, where $G = Aut\,X$\,, we can find a vector field $\tilde
{\xi}$\, on ${\bf P}^N$\, which restricts to $\xi $\, on
$X_{reg}$\,. There exists a Hermitian metric $\tilde {h}$\, on
$[H]$\,, the hyperplane line bundle on ${\bf P}^N$\,, such that
$Ric\,\tilde {h} = \tilde {\omega }$\,. (Here we denote $(Ric\,\tilde
{h})_{|U} = - \sqrt{-1}\,\partial \bar {\partial} \log (\tilde
{h}_{U})\,$ for some local expression $\tilde {h}_U$\, of the metric
$\tilde {h}$\,.) 
Hence, on $X_{reg}$\,, 
$$Ric\,\omega - \omega \, = \, Ric\, \omega - {1\over km}\,\phi^*_{L^m}\,Ric\,\tilde{h} \, = \, Ric\, \omega - Ric \,h,$$
where $h = (\tilde{h})^{1\over km}$\, is a Hermitian metric on the anticanonical bundle of $X_{reg}$\,. So, on $X_{reg}$\, we have the following equalities:
$$Ric\,\omega - Ric\,h \, = \, \sqrt{-1}\,\partial \bar {\partial}\,\log {h\over \omega ^n} \, = \, \sqrt{-1}\,\partial \bar {\partial }\, f_{\omega }$$
$$\xi (f_{\omega }) \, = \, \xi (\log {h\over \omega ^n}) \, = \, - {L_{\xi }(\omega ^n)\over \omega ^n} + h^{-1}\,L_{\xi }(h).$$
It is easy to check that
$$L_{\xi }\,(\omega ^n) \, = \, L_{\xi }\,({1\over km}\,\phi^*_{L^m}\,\tilde {\omega }^n) \, = \, \phi^*_{L^m}\,L_{\tilde {\xi }}\,({\tilde {\omega }\over km})^n.$$
Hence, on $X_{reg}$\,,
$$L_{\xi }\,(\omega ^n) \, = \, \phi ^*_{L^m}\,d\,i_{\tilde {\xi}}\,({\tilde {\omega}\over km})^n \, = \, d\,(\phi^*_{L^m}\,i_{\tilde {\xi }}({\tilde {\omega }\over km})^n).$$
We apply the Stokes' theorem and get
$$F(\xi ) \, = \,{1\over (2\,\pi)^n}\, \int_{X_{reg}} h^{-1}\,L_{\xi }\,h\,(Ric\,h)^n.$$

On the other hand, let $P' \rightarrow X$\, be the principal ${\bf
  C}^*$-bundle corresponding to $L$\,, and $\theta '$\, be the
connection (1,0)-form on $P'_{|X_{reg}}$\, corresponding to the metric
$h' = (\phi ^*_{L^m}\,\tilde {h})^{1\over m}$\,. Let further, $e_U$\,
  be a local frame of $L$\, over an open subset $U \subset X_{reg}$\,,
  and denote by $\phi _U : U \rightarrow P'_{|X_{reg}}$\, the
  corresponding local section of $P'$\,. Then, on $U$\,, which we
  identify with its image in $P'$\,under $\phi _U$\,, 
$$h^{-1}\,L_{\xi }\,(h) \, = \, {1\over k}\,(h')^{-1}\,L_{\xi }\,(h') \, = \, {1\over k}\, \theta'_{|U} (\,\xi '),$$
where $\xi '$\, is the vector field on $P'$\, (over $X_{reg}$\,) induced by $\xi $\,. Hence
$$F_U(\xi ) \, := \, {1 \over (2\pi)^n}\,\cdot {1\over
k^{n+1}}\,\int_U(h')^{-1}\,L_{\xi }\,(h')\,.(Ric\,h')^n \,$$
$$ =
\, {1\over (2\pi)^n}\cdot {1\over k^{n+1}}\, \int_U \theta '_{|U}(\xi ')\, (Ric\, h')^n$$
$$=\, {1\over n+1}\cdot {1\over (2\pi)^{n+1}}\cdot {1\over
  k^{n+1}}\,\int_U (\theta '_{|U}(\xi ') + \sqrt{-1}\Theta '_{|U}
)^{n+1},$$
where $\Theta '$\, is the curvature of $\theta '$.

Suppose now that $\nu : \tilde {X} \rightarrow X$\, is a $G$-equivariant resolution of singularities of $X$\,. Then $\xi $\, has a unique lift $\tilde {\xi }$\, to a tangent vector field on $\tilde {X}$\,, and
$$\int_U(\theta '_{|U}(\xi ') + Ric\,h')^{n+1} \, = \, \int_{\tilde U}
(\tilde {\theta '}_{|\tilde {U}} (\tilde {\xi}') + Ric\, \tilde {h'})^{n+1},$$
where $\tilde {h'} = \nu ^*(h')$\,, and $\tilde {\theta '} = \nu
^*(\theta ')$\, are the pull-backs of the corresponding objects to
$\tilde {L} := \nu ^*\,L$\, and $\tilde P' = \nu ^*P'$\,
respectively. Thus we
 get the desired formulas for the generalized Futaki invariant:
$$F(\xi ) \, = \, ({1\over 2\pi})^n \cdot ({1\over k})^{n+1}\,\int_{\tilde {X}}(\tilde {h'})^{-1}\,L_{\tilde {\xi }}(\tilde {h'})\,(Ric\,\tilde {h'})^n$$
$$=\, {1\over (n+1).k^{n+1}}\,\int_{\tilde {X}}(\tilde {\theta '}(\tilde
{\xi}') +{\sqrt{-1}\over 2\pi}\,\tilde {\Theta}'_{|\tilde {U}})^{n+1}.$$
Here ${1\over 2\pi}\, Ric\,\tilde{h'} = {\sqrt{-1}\over
  2\pi}\tilde{\Theta'}_{|\tilde {U}}$\,, where $\tilde{\Theta'}$\, is the curvature
of the metric $\tilde {h'}$.

One can apply the same argument as in \cite{Mab1} to the first of the
above formulas, and to show that the generalized Futaki invariant
vanishes on each nilpotent element of the algebra of admissible vector
fields on $X$ - the argument from the cited paper is valid also for line
 bundles which give a birational morphism of a manifold to a projective space.
 We will use the second formula in the next section.

$ $\par  

\section{``Futaki-Morita'' for non-smooth varieties}

The purpose of this section is to show that the generalized Futaki invariant, 
although analytically defined, has purely algebraic description. This invariant appears as a special case of the Futaki-Morita polynomials (cf. \cite{F-M}) 
defined as soon as we have a $H$-equivariant principal $G$-bundle over a 
(normal) variety (we work with affine algebraic groups only). Although we are
 mainly interested in the special case of the Futaki invariant, it is 
instructive to consider the general case of such polynomials.

$ $\par

\subsection{Futaki-Morita polynomials}

 Let $G$ \,be an algebraic group with Lie algebra $\g$\,. Let $I^*(G)$\, denote the subalgebra ${\bf C}[\g ^*]^G$\, of Ad$G$ -invariant polynomials over $\g$\,. The Weil algebra morphism $W_G$\, is a map
$$W_G : I^*(G) \longrightarrow H^{2*}(BG, {\bf C}),$$
where $BG$\, is the classifying space (or its algebraic approximation) of the 
principal $G$-bundles. Given a principal $G$-bundle $P \rightarrow X$\,, there is a map $\langle P \rangle : X \rightarrow BG$\, which defines a map in the cohomologies
$$\langle P \rangle ^* : H^*(BG, {\bf Z}) \longrightarrow H^*(X, {\bf Z}).$$
The composition map $w_G(P) := \langle P \rangle ^*\circ \, W_G$\, is the Weil homomorphism 
corresponding to the principal bundle $P \rightarrow X$\,.

$ $\par

Suppose we are given a principal $G$-bundle $p : P \rightarrow X$\, such that 
an algebraic group $H$\, acts on $P$\, and $X$\, commuting with the action of 
$G$\, on $P$\, and making the projection map $p$\, $H$-equivariant. By using 
the universal $H$-bundle (or an appropriate approximation) $EH \rightarrow BH$\, we can construct the principal $G$-bundle
$$p_H : P_H := EH\times _HP \longrightarrow EH\times_HX,$$
and define the corresponding map
$$\langle P_H \rangle ^* : H^*(BG, {\bf Z}) \longrightarrow H^*(EH\times_HX, {\bf Z}).$$

If $X$\, is proper (of dimension $n$), we can compose $\langle P_H \rangle ^*$\, with the 
Gyzin map
$$\gamma ^{X,H}_* : H^*(EH\times _HX, {\bf Z}) \longrightarrow H^{*-2n}(BH, {\bf Z}),$$
corresponding to the $H$-equivariant map $X \rightarrow pt$\,.

$ $\par

Suppose $X$\, is smooth $n$-manifold. In this case Futaki and Morita define a map
$$F := F_P : I^*(G) \longrightarrow I^{*-n}(H)$$
as follows (cf. \cite{F-M}). For an element $\eta \, \in \, Lie\,H$\, and $\phi \, \in \, I^{n+k}(G)$
$$F(\phi )(\eta ) :=  \int_X \phi (\Theta + \theta (\eta _*)) ,$$
where $\theta $\, is a $(1,0)$-type\,$G$-connection on $P$\,, $\Theta $\, is 
its curvature, and $\eta _*$\, is the vector field on $P$\, generated by 
$\eta $. Futaki and Morita \cite{F-M} prove that $F$\, does not depend on the choice of the connection $\theta $\, , and fits in the following commutative diagram
$$\xymatrix{
I^*(G) \ar[d]_{W_G} \ar[rr]^F && I^{*-n}(H) \ar[d]^{W_H}    \\
H^{2*}_G(pt, {\bf C}) \ar[r]^{\langle P_H \rangle ^*} & H^{2*}_H(X, {\bf C}) \ar[r]^{\gamma_*^{X,H}} & H^{2(*-n)}_H(pt, {\bf C})        }$$

Now we will show that the map $F$ can be defined, and such a diagram can be
 drawn, in the case of nonsmooth varieties too. In what follows we consider 
only $H$-equivariant principal $G$-bundles, such that the actions of $H$\, and $G$\, on them commute.

$ $\par
Notice first that for any birational morphism $\phi : Y \rightarrow Z$\,of 
proper manifolds, and principal $G$-bundle $Q \rightarrow Z$\, the morphisms 
$F_Q$\, and $F_{\phi ^*Q}$\, coincide. On the other hand, for any two $H$-
equivariant resolutions of singularities  $ \nu _i : \tilde{X}_i \rightarrow X,\,\, i=1, 2,\,$there exists a $H$-resolution $\nu : \tilde{X} \rightarrow X$\,, and maps $\mu _i : \tilde{X} \rightarrow \tilde{X}_i$\,, such that $\nu _i \circ \mu_i = \nu ,\,\, i=1, 2$.\, It follows that $F_{\nu _1^*P} = F_{\nu _2^*P} = F_{\nu ^*P}$\,, and we can define the map $F_P$\, by using any $H$-equivariant resolution of the singularities of $X$ (\cite {B-M}).

$ $\par
Let $\nu :\tilde {X} \rightarrow X$\, be a $H$-equivariant resolution of 
singularities. We then have the commutative diagram
$$\xymatrix{
H_G^{2*}(pt, {\bf C}) \ar[r]^{\langle P_H \rangle ^*} \ar [dr]_{\langle \tilde{P}_H \rangle ^*}
   &H_H^{2*}(X, {\bf C}) \ar[r]^{\gamma_*^{X,H}} \ar[d]^{\nu ^*}
   &H_H^{2(*-n)}(pt, {\bf C}),  \\
   &H_H^{2*}(\tilde{X}, {\bf C}) \ar[ur]_{\gamma_*^{\tilde{X},H}}
}$$
where $\tilde{P} := \nu ^* P$.

By the Futaki-Morita theorem we have that
$$W_H\circ F_{\tilde{P}} \, = \, \gamma _*^{\tilde{X},H}\circ \langle \tilde{P}_H \rangle ^*\circ W_G,$$
and by the diagram above we conclude that
$$W_H\circ F_P \, = \, \gamma _*^{X,H}\circ \langle P_H \rangle ^*\circ W_G.$$
It will be more convenient for our purposes to consider the map
$\tilde F := (\sqrt{-1}/ 2\pi)^n\cdot F$\, instead of
$F$. Futaki-Morita theorem can be rephrased for $\tilde F$\, as
follows:
$$\gamma_*^{X,H}\circ \langle P_H\rangle ^*\circ((\sqrt{-1}/
  2\pi)^*\, W_G) = ((\sqrt{-1}/ 2\pi)^{*-n}\, W_H)\circ
\tilde F.$$
In the special case of an almost Fano variety $X$ and a principal
bundle $\tilde P'$ as in the preceding section, we have $H :=
Aut\,X$\,, and $I^*(G) =
I^*({\bf C}^*) = {\bf C}[t]$\,. The generalized Futaki invariant of
$X$ in this case is given by $\tilde F ({1\over n+1}\cdot ({t\over k})^{n+1})$.

$ $\par
{\bf Remark.} It is well known that for a reductive group $G$\, with a
maximal torus $T$\, and Weyl group $W$\, there is an isomorphism $
I^*(G) \rightarrow I^*(T)^{W}$\,, induced by the natural restriction
map. On the other hand, for $S_{\bf Z}(T) := {\bf Sym}[\chi (T)]$\, we have
the inclusion
$$S_{\bf Z}(T)^W \subset S_{\bf Z}(T)^W\otimes {\bf C} = I^*(T)^W,$$
and the coefficient before $W_G$ \, is chosen so that, for $G =
GL(n,{\bf C})$\,, the homogeneous generators $\sigma_1,\dots ,
\sigma_n$\, of $S_{\bf Z}(T)$\, as an algebra map to the Chern classes of
the principal $G$-bundle $P_H \rightarrow M_H$, i.e. to the
$H$-equivariant Chern classes of $P$\,.

$ $\par

\subsection{Universality of the Weil morphism}

 Let $G = G_1\ltimes Rad_uG$\, be a Levi decomposition of the algebraic group $G$\,. Since $Rad_uG$\, is isomorphic to an affine space, any principal $G$-bundle can be reduced to a principal $G_1$-bundle (in $C^{\infty}$ category). In 
other words, there exist a principal $G_1$-bundle $q : Q \rightarrow X$\,, and an injection of principal bundles $i : Q \hookrightarrow P$\,, such that the 
corresponding homomorphism $i : G_1 \rightarrow G$\, coincides with the above 
embedding of $G_1$\, into $G$\,.

By using a connection on $P$\, induced by a connection on $Q$\,, one can show 
that the Weil morphism $W_G$\, factors through the algebra $I^*(G_1)$
$$\xymatrix{
I^*(G) \ar[dr]_{(\sqrt{-1}/ 2\pi)^*\cdot W_G} \ar[rr]^{res} &&
I^*(G_1) \ar[dl]^{(\sqrt{-1}/ 2\pi)^*\cdot W_{G_1}}, \\
& H_G^{2*}(pt, {\bf C})  }$$
where $res$\, is the restriction map corresponding to the embedding $i_* : \g_1 \hookrightarrow \g$\,.

The general theory of equivariant cohomology and Chow groups says that
$I^*(G_1)$\, can be identified (over ${\bf C}$)  both with
$A^*_{G_1}(pt) \cong A^*_G(pt)$\, and $H^{2*}_{G_1}(pt) \cong
H^{2*}_G(pt)$\,. After identifying $I^*(G_1)$\, with $A^*_G(pt)_{{\bf
    C}} := A^*_G(pt)\otimes _{{\bf Z}}{\bf C}$\,, denote by $i_g$\,
the resulting map from $I^*(G)$\, to $A^*_G(pt)_{{\bf C}}$\,.

$ $\par
{\bf Claim.} We have the commutative diagram

$$\xymatrix{
I^*(G) \ar[dr]_{(\sqrt{-1}/ 2\pi)^*\cdot W_G} \ar[rr]^{i_g} && A^*_G(pt)_{{\bf C}} \ar[dl]^{cl_G}.\\
   & H_G^{2*}(pt,{\bf C}) }$$

$ $\par

Proof of the claim. Without a loss of generality we may assume that
$G$ \,is reductive itself. Let $T$\, be its maximal torus, and $W =
N(T)/T$\, be the corresponding Weyl group. We start with the
commutative diagram 
$$\xymatrix{
           & I^*(G) \ar[dl]_{i_G} \ar[r]^{res^T_G} & I^*(T) \ar[rd]^{i_T} \\
A^*_G(pt)_{{\bf Q}} \ar[rd]^{cl_G} \ar[rrr]^{\phi ^*_{T,G}} &&& A^*_T(pt)_{{\bf Q}} \ar[dl]_{cl_T},   \\
&H_G^{2*}(pt, {\bf Q}) \ar[r]^{\phi ^*_{T,G}} & H^{2*}_T(pt, {\bf Q})  }$$
the horizontal arrows being injective with images the $W$-invariant subalgebras in the targets.

It is a straightforward calculation to check that $cl_T\circ i_T =
(\sqrt{-1}/ 2\pi)^*\cdot W_T$.\, On 
the other hand, we will prove in a moment that
$$\phi^*_{T,G}\circ W_G \, = \, W_T\circ res_G^T, \eqno(*)$$
which gives  $(\sqrt{-1}/ 2\pi)^*\cdot W_G = cl_G\circ i_G$\, ,
and we will be done.

To prove (*), one may use the algebraic description of the Weil morphism due to Beilinson and Kazhdan. To make the exposition self-contained, we give here the idea of the proof, following the presentation of Beilinson-Kazhdan construction as given in \cite{Esna}.

$ $\par

Recall the Atiyah extension $\A (P)$\, associated to a principal $G$-bundle $P$\, over a smooth manifold $Y$\, (cf. \cite{Atiy}):
$$\A(P) \, : \, 0 \longrightarrow Ad_Y\g \longrightarrow TP/G \longrightarrow TY \longrightarrow 0,$$
where $ Ad_Y\g := P\times _{\rho }\g$\, ($\rho $\, is the adjoint 
representation of $G$\, in $\g$\,), and $TP/G$\, is considered as a vector 
bundle on $Y$.

In their construction Beilinson and Kazhdan consider the iterated $n$-extension associated to $\A(P)^*$\, - the dual Atiyah extension :
$$\A(P)_n^* \,\, : \,\, 0 \longrightarrow \Omega _Y^n \longrightarrow \wedge^n(T^*P/G) \longrightarrow \wedge^{n-1}(T^*P/G)\otimes Ad_Y\g^* \longrightarrow \dots$$
$$\dots \longrightarrow \wedge^{n-i}(T^*P/G)\otimes S^iAd_Y\g^* \longrightarrow \dots \longrightarrow S^nAd_Y\g^* \longrightarrow 0,$$
and show that the connecting morphism
$$H^0(Y, S^nAd_Y\g^* ) \, \longrightarrow \, H^n(Y, \Omega ^n_Y),$$
when evaluated on $EG = (G^{l+1})_l \rightarrow (G^{l+1}/G)_l = BG$\, - the universal principal $G$-bundle, gives the Weil morphism
$$H^0(BG, S^nAd_{BG}\g^* ) = (S^n\g^*)^G \, \longrightarrow \, H^n(BG, \Omega ^n_{BG}) = H^{2n}_{DR}(BG).$$
As usual, for the RHS equality one needs the reductivity of $G$.

$ $\par

Suppose we are given a principal $H$-bundle $Q \rightarrow X$\,, where $H$\, is a closed subgroup of $G$\,. We have the morphism of principal bundles on $X$ : \, $Q \rightarrow Q\times _HG$\,. By the functorial properties of the Atiyah 
extension \cite{Atiy} we have the morphism of the corresponding extensions \,$\A(Q\times _HG)^* \rightarrow \A(Q)^*$\,, which in turn determines a morphism of the iterated $n$-extensions associated. Finally we have the commutative 
diagram
$$\xymatrix{
H^0(X, S^nAd_X\g^* ) \ar[rr] \ar[dr] &&H^0(X, S^nAd_X\h^*) \ar[dl] \\
&H^n(X, \Omega^n_X)         }$$
where the horizontal morphism is induced by the canonical map $\g^* \rightarrow \h^*$\,.

$ $\par

Suppose finally that $X = BH = (H^{l+1}/H)_l$\,, and $Q = EH = (H^{l+1})_l$\,. 
Since $EH\times _HG$\, is also a pull-back of $EG$\,, we have the commutative 
diagram
$$\xymatrix{
(S^n\g^*)^G \ar@{=}[r] \ar[dd]_{res^H_G} & H^0(BG, S^nAd_{BG}\g^*)
\ar[rr]^{({\sqrt{-1}\over 2\pi})^n\cdot W_G} \ar[d] && H^n(BG, \Omega^n_{BG}) \ar@{=}[r] \ar[d] &H^{2n}_{DR}(BG) \ar[dd]^{\phi _{H,G}^*} \\
& H^0(BH, S^nAd_{BH}\g^*) \ar[rr] \ar[d] && H^n(BH, \Omega ^n_{BH}) \ar[d]\\
(S^n\h^*)^H \ar@{=}[r] & H^0(BH, S^nAd_{BH}\h^*)
\ar[rr]^{({\sqrt{-1}\over 2\pi})^n\cdot W_H} && H^n(BH, \Omega ^n_{BH}) \ar@{=}[r] & H^{2n}_{DR}(BH).}$$
The right-hand-side equalities hold when $G$ \, and $H$\, are
reductive. Thus we get (*).

$ $\par

Another way of proving (*) goes as follows.

Let $EG \rightarrow BG$\, be the universal $G$-bundle, and let $B$
\,be a Borel subgroup of $G$\,. Consider the corresponding fibration
$$\varphi : EG/B \rightarrow BG\,$$
 with fibre  $G/B\,$. Then, 
$$\varphi ^*(EG) = EG\times _BG\,,$$
 where $B$\, acts on $G$\, from the left, and the
cohomology morphism $$\varphi ^* : H^*(BG, {\bf C}) \rightarrow H^*(EG/B,
{\bf C})\,$$
 coincides with the morphism $$\phi ^*_{B,G} : H^*_G(pt,
{\bf C}) \rightarrow H^*_B(pt, {\bf C})\,.$$

By the construction, the structure group of the principal $G$-bundle
$\varphi ^*(EG)$\, is reducible to the subgroup $B$\,, and since the
map $\,\,\,\varphi ^*\circ W_G = \omega _G(\varphi ^*(EG))$\,\,\, does not
depend on the connection chosen, we take a $B$-connection on $EG
\rightarrow EG/B$\, and induce a $G$-connection on $\varphi ^*(EG) \rightarrow
EG/B$\,. This induced connection gives $$\omega _G(\varphi ^*(EG)) = W_B\circ
res^B_G\,,\,$$ or in the other notations $$\phi ^*_{B,G}\circ W_G =
W_B\circ res^B_G\,. \eqno(1)$$

\par

We may apply the arguments above also to the bundle $EG \rightarrow
EG/B$\, and to the maximal subtorus $T$\, of $B$\,. The principal
$B$-bundle $EG$\, is reducible to a $T$-principal subbundle because
the manifold $B/T$\, is isomorphic to an affine space, and therefore
the fibration $EG/T \rightarrow EG/B$\, has a section. Thus we get 
$$W_T\circ res^TB = \phi ^*_{T,B}\circ W_B. \eqno(2)$$
Now (*) easily follows from (1) and (2).

$ $\par

\subsection{Inserting a row in the Futaki-Morita diagram}

We can apply the same reasoning as in the case of equivariant cohomology groups above to a $H$-equivariant principal $G$-bundle $p : P \rightarrow X$\,, with commuting actions of $H$\, and $G$\, on $P$. We obtain a map 
$$[P_H]^* : A^*_G(pt)_{{\bf Q}} \longrightarrow A^*_H(X)_{{\bf Q}},$$
where $[P_H] : X_H \rightarrow BG$\, is the classifying map corresponding to 
the principal $G$-bundle $p_H : P_H \rightarrow X_H$\,. When $X$\, is smooth 
and proper, we can combine the Futaki-Morita polynomials and the equivariant Chow groups in the following diagram:
$$\xymatrix{
I^*(G) \ar[rrrr]^{\tilde F} \ar[rd]^{i_G} \ar[dd]_{({\sqrt{-1}\over
    2\pi})^*\cdot W_G} &&&& I^{*-n}(H) \ar[dl]_{i_H}
\ar[dd]^{({\sqrt{-1}\over 2\pi})^{*-n}\cdot W_H} \\
& A^*_G(pt)_{{\bf C}} \ar[dl]_{cl_G}^{\cong} \ar[r]^{[P_H]^*} & A^*_H(X)_{{\bf C}} \ar[r]^{\int_X^H} \ar[d]^{cl_H^X} & A^{*-n}_H(pt)_{{\bf C}} \ar[dr]^{cl_H}_{\cong} \\
H_G^{2*}(pt, {\bf C}) \ar[rr]^{\langle P_H \rangle ^*} && H_H^{2*}(X, {\bf C}) \ar[rr]^{\gamma_*^{X,H}} && H_H^{2(*-n)}(pt, {\bf C}).  }$$
When $X$ is a normal algebraic variety we have the same diagram but without the map  $cl_H^X$\,. This is easily seen by using $H$-desingularizations. (The map $\int _X^H\circ [P_H]^*$\, does not depend on the desingularizations and, under the identifications in the diagram, coincides with $\gamma _*^{X,H}\circ \langle P_H \rangle ^*$\,.)

$ $\par

Choose a Levi decomposition $H = H_1\ltimes Rad_u\,H$\, of the group
$H$\,. According to this decomposition we have
$$I^*(H) \subset {\bf C}[\h _1]\otimes {\bf C}[rad_u\h] = {\bf C}[\h
_1][rad_u \h].$$
By its definition, $i_H$\, is given by sending $\varphi \in I^*(H)$\,
to its degree 0 component in the latter algebra
$$\varphi _0 := i_H(\varphi ) \in {\bf C}[\h _1].$$
Since $\varphi $\, is $Ad\,H$-invariant, and since all Levi subgroups
of $H$\, are conjugate, the property that the difference $\varphi -
\varphi _0$\, is $0$\, is well defined : if it holds for some Levi
subgroup $H_1$\, of $H$\,, then it holds for all such subgroups.

\begin{theorem}\label{algebranature}
Let $G$\, and $H$\, be algebraic groups, $p : P \rightarrow X$\, a $H$-
equivariant principal $G$-bundle over a normal variety $X$\,. Let the actions 
of $G$\, and $H$\, on $P$\, commute. Given $\phi \, \in \, I^*(G)$\,,
denote by $\alpha $ its image in $A^*_H(X)_{{\bf C}}\,\, : \,\,  \alpha
= [P_H]^*\circ i_G(\phi)\,$. Then,

1. The polynomial $\tilde F(\phi )$\, is {\bf 
algebraic}, i.e. for each $\eta \in \h = Lie\,H$ \,\, $\tilde F(\phi )(\eta) =
f_{\alpha }( \eta )\, $, iff $\tilde F(\phi ) - \tilde F(\varphi )_0 =
0$\, for some Levi decomposition $H = H_1\ltimes Rad_uH$\, of $H$ ;

2. In particular, if $H$ \,is reductive, then $\tilde F(\phi )$\, is algebraic for all $\phi \, \in \, I^*(G)$\,.

\end{theorem} 

Note that for degree 1 elements  $\varphi \in I^*(H)$\, the property
$\varphi - \varphi _0 = 0$\, means exactly that $\varphi $\, vanishes
on the nilpotent elements of $\h$\,.As we have already seen in the
previous section, the generalized Futaki  
invariant of an almost Fano variety $X$\, vanishes on all nilpotent
elements of $Lie\, H$\,, where $H = Aut(X)$\,. Combining this with the
Remark in the 
subsection 3.1 it follows that

\begin{coroll}
The generalized Futaki invariant is algebraic. It can be represented
as ${1\over n+1}\,f_{\alpha }$\,, where $\alpha = (c^H_1(X))^{n+1}$\,.
\end{coroll}

{\bf Remarks.}1. According to the remark in the end of Section 2, we
have that the generalized Futaki invariant can be lifted to a group
character. Another approach to see this is by using the Chern-Simons
invariants. A deep study of this approach can be found in \cite{F-M}.

 2. In the notations above let further $T$\, be a
maximal torus of $H$\, which lies in $H_1$\,, and $T = T'\times
T''$\,, where $T'$\, is a maximal torus in $H_1'$\, and $T'' =
Z(H_1)_0$
 - the connected component of the unit of the center of $H_1$\,. Let
$W = W(H_1',T')$\, be the corresponding Weyl group. By the structure
theorem for the Chow cohomology groups we have
$$A^*_H(pt)_{{\bf C}} \cong A^*_T(pt)_{{\bf C}}^W = I^*(T)^W.$$
On the other hand, a theorem of Chevalley says that the restriction
map $I^*(H_1) \rightarrow I^*(T)^W$\, is an isomorphism of graded
algebras. Hence, each element $\alpha \in A^*_H(pt)_{{\bf Q}}\,$
defines, implicitly, a polynomial from $I^*(H_1)$\, and, by using the
canonical morphism $\h \rightarrow \h _1$\,, an element of $I^*(H)$\,.

3. Consider the subalgebra $\z(\h _1) = \t ''$\, of $\h $\,. The
   corresponding restriction morphism
$$res^{\h }_{\t ''} : I^*(H) \rightarrow I^*(T'') = {\bf C}[\t '']$$
provides a correct definition of the property $\,\,\varphi _1 := \varphi -
res^{\h }_{\t ''}\,\varphi $\,\,\, being zero. In the particular case of
$\varphi \in I^1(H)$\,, this property is equivalent to $\varphi $
\,being a character of $\h $\, which vanishes on its nilpotent
elements. Finally, the canonical projection $\,\,\h \rightarrow \z (\h _1)
= \t ''$\,\,\, determines an embedding
$${\bf C}[\t ''] \hookrightarrow I^*(H)\,,$$
which gives a way of finding explicitly the polynomials in $I^*(H)$\,
corresponding to the elements $f_{\alpha } \in {\bf Sym}[\chi (T'')]_{{\bf
    Q}}\,$. These polynomials are $f_{\alpha }$\, themselves.

\section{Examples.}

We give in this section examples of calculating the generalized Futaki
invariant $F$\, of some almost Fano varieties $X$. Since $F$\, is an
equivariant Chow class over a point obtained by integration over
$X$\,, the main tool we use, when it cannot be calculated directly, is
the Bott-type residue formula proved by Edidin and Graham \cite
{E-G2}. This formula can also be used for calculating the polynomials
$\tilde F (\phi )\,, \phi \in I^{n+k}(G)$\,, provided that they are
algebraic in the sense of Theorem \ref{algebranature}.

Suppose $\tilde F (\phi )$\, is algebraic. As we have already proved,
to see the behaviour of $\tilde F (\phi )$\, on the elements of $\h$\,
it is enough to  fix a maximal torus $T$\, of $H$\, and find the
element of $A^k_T(pt)$\, corresponding to $\tilde F (\phi )$\,. As we
have shown, when $k=1$\, this element defines an element of $\h ^*$\,
representing $\tilde F (\phi )$\, as a functional on $\h$\,.

If one has a good description of a $H$-desingularization $\tilde X$\,
of $X$\, (as it is in the case of toric varieties), he could calculate
$\tilde F (\phi )$\, on $\tilde X$\,. There are some cases though
where, by using the Bott-type residue formula, it is easier to
calculate $\tilde F (\phi )$\, on $X$\, itself. For one such example
with $\tilde F (\phi ) = F$\, see below. 

{\bf Remark.}  In the general case of a smooth K\"ahler manifold $M$\,
and a polynomial $\tilde F (\phi )$\, there are formulas proved by
Futaki-Morita \cite{F-M} \,  Futaki \cite{Fut2}, and Tian \cite{Tia2}  
which calculate
the value of $\tilde F (\phi )$\, on nondegenerate  holomorphic vector
fields $\xi $\, on $M$.\, When $\xi $ lies in a Lie subalgebra $\k$\,
corresponding to a compact subgroup $K$\, of $H$\,, these are
precisely the Bott-type residue formulas we use. In fact, one can
prove that for any compact subgroup $K \subset H$\, the restriction of
$\tilde F (\phi )$ \, to $Lie\,K$\, can be represented by a
$K$-equivariant cohomology class of a point.

$ $\par

We recall now the essence of the Bott residue formula for algebraic
varieties as presented in \cite{E-G2}. 

$ $\par

Let $T$\, be an algebraic torus which acts on an algebraic scheme
$X$\,, and let $M$ \,be a smooth algebraic $T$-manifold. Denote by
$X^T$\, (resp. $M^T$\,) the set of connected components of the fixed
points on $X$\, (resp. $M$\,) under the action of $T$\,. Suppose 
$$f : X \longrightarrow M$$
is a $T$-equivariant embedding such that $X^T \subset M^T$\, (meaning
that if $F \in M^T$\, and $F \cap  X \neq \varnothing$\,, then $F \in
X^T$\,). For each $F \in X^T$\, denote by $i_F : F \rightarrow X$\,
(resp. $j_F : F \rightarrow M$\.) its embedding into $X$\,
(resp. $M$)\,. Hence, for each such $F$\, we have $j_f = f\circ
i_F$\,. Denote finally by ${\Q}$\, the quotient ring of ${\bf
  Sym}[\chi (T)]$\,. Then, Proposition 6 of \cite{E-G2} says that the
map  
$$f_* : A_*(X)\otimes {\Q} \longrightarrow A_*^T(M)\otimes {\Q}$$
is an inclusion, and for $\alpha \in A_*^T(X)\otimes {\Q}$\, we
have
$$\alpha = \sum_{F\in X^T}i_{F*}\,{j^*_F\circ f_*{\alpha}\over
  c^T_{top}(N_FM)}\,\,,$$
where $N_FM$ \, denotes the normal sheaf of $F$\, in $M$\,, and 
$c_{top}^T(N_FM)$\, denotes the top Chern $T$-equivariant class of that
sheaf. Recall that by a theorem of Iversen \cite{Iver}, the components
of the fixed-point set of a torus action on a smooth algebraic variety
are smooth submanifolds.

In particular,
$$[X]_T = \sum_{F \in X^T}i_{F*} \beta _F,$$
where
$$\beta _F = {j_F^*\circ f_*[X]_T\over c^T_{top}(N_FM)}\,\,.$$
The algebraic Futaki-Morita polynomials $\tilde F (\phi )$\, of $X$\,
can be expressed as follows
$$\tilde F (\phi ) = \int_X^T\, p^T(P) = \pi _{X*}^T\,(p^T(P) \cap
[X]_T)$$
$$= \sum_{F \in X^T}\pi _{F*}\,(p^T(i_F^*\,P)\,\cap \beta _F )\,,$$
and to find $\tilde F (\phi )$\, one needs only to know $\beta _F$\,
for each $F \in X^T$\,. Here $p^T(P) := [P]^\circ i_G(\phi )\,$ . When 
$ G = GL(r, {\bf C})$\, $p^T(P)$\, is just a polynomial of the Chern
classes of the principal bundle $P$\,.

In the special case when $i_F$\, and $f$\, are regular embeddings,
e.g. when $F \cap Sing\,X = \varnothing$\,, we have
$$j_F^*\circ f_*\,[X]_T = i_F^*\circ f^*\circ f_*\,[X]_T =
i_F^*(c^T_{top}(N_XM)\cap [X]_T)$$
$$= c^T_{top}(i_F^*\,N_FM)\cap [F]_T\,,$$
as well as
$$c_{top}^T(N_FM) = c^T_{top}(N_FX) \cup c_{top}^T(i_F^*N_FM)\,.$$
Hence in this case
$$\beta _F = {[F]_T\over c_{top}^T(N_FX)}\,.$$

Another important case in which one has an effective method of
computing $\beta _F$\, is when $F$\, is an attractive isolated fixed
point. Although it happens that all fixed points in our examples are
attractive, we will not use this method, referring to \cite{Bri1} and
\cite{Bri2}  for details about it.
$ $\par 

1. {\bf Complete intersections.} As a first example we prove a formula
for the generalized Futaki invariant $F$ of complete intersections in
${\bf P}^N$\,. In this case, since all ingredients of the equivariant
class representing $F$\, can be expressed in terms of equivariant
classes of vector bundles on ${\bf P}^N$\,, it can be calculated
directly  without use of residue formulas. The formula we prove has
already appeared in the recent preprint by Zhiqin Lu \cite{Lu}, where he proves
it by using technique different from ours.

$ $\par

Let $X \hookrightarrow {\bf P}^N$\, be a complete intersection given
by the polynomials $f_1,\dots , f_k$\, of degree $\deg\,f_i = d_i\,, i
= 1,\dots ,k.$\, Denote by $H := Aut(X \subset {\bf P}^N)\, <
\,Aut({\bf P})^N$\, the group of automorphisms of $X$\, in ${\bf
  P}^N$\,. Let $T \,<\, H$\, be a maximal torus of $H$\,, and suppose
it is $m$-dimensional. We choose
coordinates $(X_0 :\cdots : X_N)$\, of ${\bf P}^N$\, such that the
action of $T$ on ${\bf P}^*N$\, is given by
$$(u_1,\dots , u_m)\circ (X_0:\cdots : X_N) := (u^{\gamma_0}X_0:\cdots
:u^{\gamma_N}X_N )\,,$$
where $u := (u_1,\dots , u_m) \in T$\,, for $i=0,\dots , N$\, 
we denote $u^{\gamma _i} = \prod _{j=1}^m\,u_j^{\gamma _i^j}$\,, and
$\gamma _0 +\cdots + \gamma _N = 0$\,. 

Without a loss of generality we may assume that $f_1,\dots , f_k\,$
are $T$-semi-invariants, i.e. that $u\circ f_i = u^{k^i}\,f_i$\, for
$k^i \in {\bf Z}^m\,$.  

Let ${e_1,\dots , e_m}\,$ be a basis of $\t  = Lie\,T$\,, and
  ${e^1,\dots , e^m}\,$ be its dual basis of $\t ^*$\,. The Futaki
  invariant can be interpreted as an element of $\t ^* $\, and in the
  basis chosen it has the form
$$ F = a_1\,e^1 + \cdots + a_m\,e^m\,.$$
Since the component $a_i\,e^i\,$ is the restriction of $F$\, to the
$i$-th factor ${\bf C}\,e_i\,$ of $\t $\,, we reduce the calculation
of $F$ \, to the case of 1-dimensional torus action $T = {\bf
  C}^*$\,. So, we have now a one-dimensional torus $T$\, acting on
${\bf P}^N$\, via
$$u\circ (X_0:\cdots : X_N ) := ( u^{\gamma ^i_0}\,X_0 :\cdots :
u^{\gamma ^i_N}\,X_N)\,,$$
for some integer $\gamma ^i_j\,\,\,, j = 1,\dots , N\,$\, and
$$u\circ f_j = u^{k_i^j}\,f_j\,\,\,\,\,\,\,\,\,\,\,\,\,\,\,\,\,\,\, j
= 1,\dots , k\,.$$

It is well known (see for example \cite{E-G1}) that the $T$-equivariant
Chow ring of the projective space under chosen action is given by
$$A_T^*({\bf P}^N) = {\bf Z}[h, t]/{\prod _{j=0}^N(h + \gamma
  _j^i\,e^i)}\,,$$
and that
$$[X]_T = \prod _{j=1}^k\,(d_j\,h + k_i^j\,e^i) \,\,\,\,\in \,\,
A_T^k({\bf P}^N)\,.$$
On the other hand, since $\gamma _0^i + \cdots + \gamma _N^i = 0\,$
the equivariant Chern class $c_1^T({\bf P}^N) = (N+1)\,h$\,, and by
the adjunction formula we have
$$c_1^T(X) = c_1^T({\bf P}^N) - c^T_k(N_XM)$$
$$= (N+1)\,h - (d\,h + k_i\,e^i)\, ,$$
where $d = d_1 + \cdots + d_k$\,, and $k_i = k_i^1 + \cdots +
k_i^k$\,.

The generalized Futaki invariant is given now by
$$F = {1\over N-k+1}\, \int_X^T\,((N+1-d)\,h - k_i\,e^i)^{N+1-k}$$
$$= {1\over N-k+1}\, \left[((N-d+1)\,h - k_i\,e^i)^{N-k+1}\cdot \prod
_{j=1}^k(d_j\,h + k_i^j\,e^i)\right]_{h^N}\,.$$
Thus we get that in the basis $(e^1,\dots , e^m)$\, the $i$-th
coordinate of $F$\, is given by the formula
$$a_i = (N-d+1)^{N-k}\cdot \prod
_{j=1}^kd_j\,\sum_{j=1}^k\,\left ( {N-d+1\over N-k+1}\cdot{1\over d_j} - 1
\right ) \,k_i^j\,.$$ 

$ $\par

2.{\bf Toric examples.} For at least three reasons the Bott-type 
formula is very  convenient for computing the generalized Futaki
invariant of toric varieties. First, all toric varieties have toric
resolutions, which have very clear combinatorial explanation. Second,
the torus action on a toric variety has only isolated fixed
points. Third, the equivariant line bundles with respect to the torus
action are well understood. 

We give below two examples of almost Fano toric varieties and compute
their Futaki invariants. The reason we have chosen these will be
explain further on.

$ $\par

2.1. Let $X$\, be the blow up of ${\bf P}^3$ in a curve given by the
ideal $(X_0\,X_2, X_0\,X_3, X_2\,X_3)$\, ($(X_0: X_1: X_2: X_3)$\, are
the coordinates of ${\bf P}^3$\,). Then $X$\, is a Gorenstein toric
almost Fano variety with two terminal singular points. Denote by $T$\,
the torus according to which $X$\, is toric. This variety
has a small resolution $\tilde X$\,, and for its Futaki invariant
$F$\, we have
$$F = {1\over 4}\, \int_{X}^T\,c_1^T(X)^4 = {1\over 4}\,\int_{\tilde X}^T(f^*c_1^T(X))^4$$
$$= {1\over 4}\,\int_{\tilde X}^Tc_1^T(\tilde X)^4 = {1\over
  4}\,\sum_{P \in \tilde X^T} i_P^*(K_{\tilde X}^*)^4\,. \beta _P \,.$$
So, we proceed further with computations on $\tilde X$\,.

A combinatorial description of $\tilde X$\,. If we start with the fan
in $N_{{\bf R}} = {\bf R}^3$\,
$$\langle e_1, e_2, e_3\rangle, \langle e_0, e_1, e_2\rangle, \langle
e_0, e_1, e_3\rangle, \langle e_0, e_2, e_3\rangle$$
for ${\bf P}^3$\,, where $e_0 + e_1 + e_2 + e_3 = 0$\,, and $e_0, e_1,
e_2, e_3$\, span ${\bf R}^3$\,, then $\tilde X$\, has the fan
$$\sigma _1 = \langle e_1, e_2, e_2+e_3\rangle ,\, \sigma _2 = \langle
e_1, e_3, e_2+e_3\rangle\,, \sigma _3 = \langle e_3, -e_1-e_2,
e_2+e_3\rangle$$
$$\sigma _4 = \langle -e_1-e_2-e_3, e_2, -e_2-e_3\rangle\,, \sigma _5
= \langle e_1, e_2, -e_2-e_3\rangle\,, \sigma _6 = \langle
-e_1-e_2-e_3, -e_1-e_2, -e_2-e_3\rangle$$ 
$$\sigma _ = \langle -e_1-e_2-e_3, e_2, -e_1-e_2\rangle \,, \sigma _8
= \langle e_2, -e_1-e_2, e_2+e_3\rangle\,, \sigma _9 = \langle e_1,
e_3, -e_2-e_3\rangle$$
$$\sigma _{10} = \langle e_3, -e_2-e_3, -e_1-e_2\rangle\,.$$ 
The anticanonical line bundle of $\tilde X$\, is given then by the
vectors in $M = N^*$\,
$$m(1) = -e^1-e^2, m(2) = -e^1-e^3, m(3) = e^1-e^3, m(4) =
-e^2+2\,e^3, m(5) = -e^1-e^2+2\,e^3,$$
$$m(6) = e^2, m(7) = 2\,e^1-e^2, m(8) = 2\,e_1-e^2, m(9) =
-e^1+2\,e^2-e^3, m(10) = -e^1+2\,e^2-e^3.$$
($m(i)$\, is determined by $\sigma _i$\,) Denote by $P_1, \dots ,
P_{10}$\, the fixed points of the action of $T$\,on $\tilde X$\,
corresponding to $\sigma _1,\dots ,\sigma _{10}\,$ respectively. An
easy computation shows that, in the notations above and if we denote
$\beta _i := \beta _{P_i}$\,,
$$\beta _1 = {-1\over e^1\,(e^2-e^3)\,e^3}\,\,,\,\, \beta _2 = {1\over
  e^1\,e^2\,(e^2-e^3)}\,\,,\,\, \beta _3 ={-1\over
  e^1\,(e^1-e^2)\,(e^1-e^2+e^3)}\,,$$
$$ \beta _4 = {1\over e^1\,(e^2-e^3)\,e^1-e^3)}\,\,,\,\, \beta _5 =
{1\over e^1\,(e^2-e^3)\,e^3}\,\,,\,\, \beta _6 = {-1\over
  (e^1-e^2)\,(e^2-e^3)\,(e^1-e^2+e^3)}\,,$$
$$\beta _7 = {1\over (e^1-e^2)\,(e^1-e^3)\,e^3}\,\,,\,\, \beta _8 =
{-1\over e^1\,e^3\,(e^1-e^2+e^3)}\,\,,\,\, \beta _9 = {-1\over
  e^1\,e^2\,(e^2-e^3)}\,,$$
$$ \beta _{10} = {1\over e^1\,(e^1-e^2)\,(e^1-e^2+e^3)}\,.$$

Now we are ready to compute the Futaki invariant of $X$:
$$F = {1\over 4}\,\sum_{i=1}^{10}m(i)^4\,. \beta _i =
4\,(-e^1-e^2+e^3)\,.$$

$ $\par

2.2. The second example of an almost Fano toric variety $X$\, is the
     blow up of ${\bf P}^3$ in a curve with ideal $(X_0\,X_1,
     X_0\,X_2, X_1\,X_2\,)$ . This variety has three terminal singular
     points. It has also a small resolution $\tilde X$ with 12 fixed
     points. Starting with the same fan for ${\bf P}^3$ we get a fan
     for $\tilde X$
$$\sigma _1 = \langle -e_1-e_2-e_3, e_1, e_1+e_2\rangle \,\,,\,\,
\sigma _2 = \langle -e_1-e_2-e_3, e_2, e_1+e_2\rangle \,\,,\,\, \sigma
_3 \langle -e_1-e_2-e_3, e_2, e_2+e_3\rangle \,\,,$$
$$\sigma _4 = \langle -e_1-e_2-e_3, e_3, e_2+e_3\rangle \,\,,\,\, \sigma
_5 = \langle -e_1-e_2-e_3, e_3, e_1+e_3\rangle \,\,,\,\, \sigma _6 =
\langle -e_1-e_2-e_3, e_1, e_1+e_3\rangle ,$$
$$\sigma _7 = \langle e_1, e_1+e_2, e_1+e_2+e_3\rangle \,\,,\,\,
\sigma _8 = \langle e_1, e_1+e_3, e_1+e_2+e_3\rangle \,\,,\,\, \sigma
_9 = \langle e_2, e_1+e_2, e_1+e_2+e_3\rangle \,\,,\,\,$$
$$\sigma _{10} = \langle e_2, e_2+e_3, e_1+e_2+e_3\rangle \,\,,\,\,
\sigma _{11} = \langle e_3, e_1+e_3, e_1+e_2+e_3\rangle \,\,,\,\,
\sigma _{12} = \langle e_3, e_2+e_3, e_1+e_2+e_3\rangle \,\,.$$

Proceeding now as in the previous example, we compute that
$$F = 4\,(e^1 + e^2 + e^3 ).$$

$ $\par

{\bf Remark.} The generalized Futaki invariant of an almost Fano toric
variety has a very nice geometric interpretation. As an element of
$M_{{\bf R}}$\, $F$\, coincides with the barycentre of the polytope
corresponding to the anticanonical sheaf of that variety, with respect
to the Duistermaat-Heckman's measure on it.(See \cite{Mab2} for the
smooth case, and \cite{Yoto} for the general case.) By using this
interpretation the author  \cite{Yoto} has computed $F$\, for the
above toric examples getting the same result.

$ $\par

3.{\bf Example.} As a last example consider the blow up $X$\, of ${\bf
  P}^3$\, in the curve given by the ideal $(X_0\,X_3, X_1\,X_3,
  X_0\,X_2 - X_1^2)\,$. Then $X$\, is an almost Fano variety with one
  singular point. We can naturally consider $X$\, as a subvariety of
  $M := {\bf P}^3\times {\bf P}^2$\,. Denote by $F$\, the
  corresponding embedding. Let $(Y_0: Y_1: Y_2)$\, be the coordinates
  of ${\bf P}^2$\,.

The group $H := Aut\,X$\, has a two dimensional maximal torus $$T :=
\{\,\, (z_1, z_2) \,\,|\,\, z_i \in {\bf C}^*\,\,\,, i=1, 2\,\, \}\,
.$$ 
If we let $T$\, act on $M$\, via
$$(z_1, z_2)\circ (X_0 ; X_1: X_2: X_3; Y_0: Y_1: Y_2)$$ $$ := (X_0:
z_1\,X_1: z_1^2\,X_2: z_2\,X_3; z_2\,Y_0: z_1z_2\,Y_1: z_1^2\,Y_2),$$
then the embedding $f$\, will be $T$-equivariant. The action of $T$\,
on $X$\, has seven fixed points. As points of $M$\, they have
coordinates
$$P_1(1: 0: 0: 0; 1: 0: 0)\,\,,\,\, P_2(1: 0: 0: 0; 0: 0: 1)\,\,,\,\,
P_3(0: 1: 0: 0; 0: 0: 1)\,$$
$$P_4(0: 0: 1: 0; 0: 1: 0)\,\,,\,\, P_5(0: 0: 1: 0; 0: 0: 1)\,\,,\,\,
P_6(0: 0: 0: 1; 1: 0: 0)\,\,$$
$$P_7(0: 0: 0: 1; 0: 1: 0)\,\,.$$
The point $P_5$\, is the only singular point of $X$\,. Denote by $i_j
: P_j \hookrightarrow X\,$ the corresponding embeddings. Since the
action of $T$\, on $M$\, has only finite number of fixed points, we
have an inclusion $X^T \subset M^T$\, and can apply the Edidin-Graham
theorem for computing $F$\,:
$$F = {1\over 4}\,\sum_{j=1}^7\,i_j^*(K_X^*)\,.\beta _j\,.$$
The calculation of $m(j) := i_j^*(K_X^*)$\, and $\beta _j$\, for
$j\neq 5$\, is straightforward.  It happens that $P_5$\, is an
attractive fixed point. So, $\beta _5$ could be calculated by using the
method from \cite{Bri1}. We here avoid finding $\beta _5$\, in
computing $F$\, proceeding in the following way. By the Bott residue 
theorem we have
$$38 = (-K_X)^3 = \sum_{j=1}^7\,m(j)^3\,.\beta _j.$$
Hence, we can exclude the term $\beta _5$\, from the formula for
$F$\,:
$$F = {1\over 4}\,\left[ \,\,38\,m(5) + \sum_{j\neq 5}\,m(j)^3\,.\,\beta
  _j\,\,(\,m(j) - m(5)\,)\,\,\right]\,.$$
Then, the term $m(5)$\, is uniquely determined by the fact that the
sum on the RHS is a linear polynomial of the generators $e^1, e^2$\,
of the character group $\chi (T)$\, of $T$\,. Notice that we don't
actually need to know $\beta _5$\, to compute $F$\,.

In this way we get the input data:
$$\beta _1 = {1\over e^1\,e_2\,(2\,e^1-e^2)}\,\,,\,\, \beta _2 =
{-1\over 2\,(e^1)^2\,(2\,e^1-e^2)}\,\,,\,\, \beta _3 = {1\over
  (e^1)^2\,(e^1-e^2)}\,\,$$
$$\beta _4 = {1\over e^1\,(e^1-e^2)\,(2\,e^1-e^2)}\,\,,\,\, \beta _6
= {-1\over e^1\,e^2\,(2\,e^1-e^2)}\,\,,\,\, \beta _7 = {-1\over
    e^1\,(e^1-e^2)\,2\,e^1-e^2)}$$

$ $\par 

$$m(1) = 3\,e^1 \,\,,\,\, m(2) = e^1+e^2\,\,,\,\, m(3) =
-e^1+e^2\,\,,\,\, m(4) = -2\,e^1\,\,$$
$$m(5) = -3\,e^1+e^2\,\,,\,\, m(6) = \,3\,e^1-2\,e^2\,\,,\,\, m(7) =
2\,(e^1-e^2)\,\,$$
and compute the output
$$F = 4\,(\,3e^1 - e^2\,).$$

$ $\par

We conclude this section by explaining why the last three examples are
interesting to the author.

Recall (\cite{Tia2}) that a special degeneration of a Fano manifold
 $M$\, is called
a morphism $\pi : W \rightarrow \Delta \,$, without multiple fibres,
 where $\Delta $\, is the one dimensional unit disc, such that

1) the relative anticanonical sheaf of $\pi $\, is ample,

2) the family of dilatons $\phi (\lambda ) : \Delta \rightarrow \Delta
   ,\,\,\,\,t \mapsto \lambda t ,\,\,\,\, 0 < |\lambda | \leq 1,\,\,$
   extends to a group $\Phi (\lambda ) \,$ of transformations of $W$,

3) the central fibre $W_0$\, of $\pi $\, is a normal almost Fano
   variety, and all other fibres are isomorphic to $M$.

\par
Since $W_0$\, is $\Phi (\lambda )$-invariant, there is defined an
admissible vector field $ v_W := -\Phi '(1)\,$ on it. A deep theorem by
Ding and Tian \cite{D-T}, strengthened by Tian \cite{Tia2}, says that 
if $M$\, is Einstein-K\"ahler, then the real part of the generalized 
Futaki invariant $F$, evaluated on $v_W$ is nonnegative.

The last three examples from this section can be realized as central
fibres of special degenerations of the Fano manifold $V_{38}$\, - the
blow-up of the projective three-dimensional space in a twisted cubic
curve \cite{Yoto}. One can show that in the all three cases the real 
parts of the corresponding evaluations are positive.

\bibliographystyle{amsalpha}

\end{document}